\documentclass[12pt,reqno]{amsart}
\usepackage{amssymb}
\usepackage{amsmath}
\usepackage{graphicx}
\usepackage[dvips,line]{xy}

\date{26 September 2007; minor revisions 19 October 2007}

\theoremstyle{plain}
\newtheorem{thm}{Theorem}[section]
\newtheorem{conj}[thm]{Conjecture}
\newtheorem{lemma}[thm]{Lemma}
\newtheorem{cor}[thm]{Corollary}
\newtheorem{quest}[thm]{Question}
\newtheorem{prp}[thm]{Proposition}
\theoremstyle{definition}
\newtheorem{exm}[thm]{Example}
\newtheorem{rmk}[thm]{Remark}

\numberwithin{equation}{section}

\newcommand{\myspread}{0.5\jot}

\newcommand\Z{\mathbb{Z}}

\newcommand\R{\mathbb{R}}
\newcommand{\mtset}{\varnothing}
\newcommand{\aff}[1]{\widetilde{#1}}
\newcommand{\rev}[1]{#1^*}
\newcommand{\trim}[1]{#1{\downarrow}}
\newcommand{\ltrim}[1]{{\downarrow}#1}
\newcommand{\vph}{\vphantom{1}}
\newcommand{\br}[1]{\langle #1\rangle}
\newcommand\eps{\varepsilon}
\newcommand{\hr}{\widetilde{\alpha}}
\newcommand{\cbox}[2]{\hbox to #1pt{\hfill$#2$\hfill}}
\newcommand{\xtal}{crystallographic}
\newcommand{\Branden}{Br\"and\'en}
\newcommand{\fake}[1]{#1^{\text{fa}}}
\renewcommand{\le}{\leqslant}
\renewcommand{\ge}{\geqslant}

\newcommand{\st}{Steinberg torus}
\newcommand{\sti}{Steinberg tori}
\newcommand{\rst}{reduced Steinberg torus}
\newcommand{\rsti}{reduced Steinberg tori}
\newcommand{\stor}[1]{\Sigma_T(\aff{#1})}
\newcommand{\storr}[2]{\Sigma_T(\aff{#1}_{#2})}
\newcommand{\rstor}[1]{\Sigma'_T(\aff{#1})}
\newcommand{\rstorr}[2]{\Sigma'_T(\aff{#1}_{#2})}

\newcommand{\des}{d}
\newcommand{\ades}{\aff{d}}
\newcommand{\Des}{D}
\newcommand{\aDes}{\aff{D}}

\DeclareMathOperator{\pk}{pk}
\DeclareMathOperator{\xpe}{epk}
\DeclareMathOperator{\lpe}{lpk}

\title{Affine descents and the Steinberg torus}

\author[K. Dilks]{Kevin Dilks}
\address{Department of Mathematics, University of Michigan,
  Ann Arbor MI 48109--1043}
\email{kdilks@umich.edu}

\author[T. K. Petersen]{T. Kyle Petersen}
\address{Department of Mathematics, University of Michigan,
  Ann Arbor MI 48109--1043}
\email{tkpeters@umich.edu}

\author[J. R. Stembridge]{John R. Stembridge}
\address{Department of Mathematics, University of Michigan,
  Ann Arbor MI 48109--1043}
\email{jrs@umich.edu}
\thanks{The work of the third author was supported in part
  by an NSA Grant.}

\begin{document}

\begin{abstract}
Let $W\ltimes L$ be an irreducible affine Weyl group
with Coxeter complex $\Sigma$, where $W$ denotes the associated
finite Weyl group and $L$ the translation subgroup.
The \st{} is the Boolean cell complex obtained by taking the
quotient of $\Sigma$ by the lattice $L$.
We show that the ordinary and flag $h$-polynomials of the \st{}
(with the empty face deleted) are generating functions over $W$ for
a descent-like statistic first studied by Cellini. We also show that
the ordinary $h$-polynomial has a nonnegative $\gamma$-vector,
and hence, symmetric and unimodal coefficients.  In the classical
cases, we also provide expansions, identities, and generating
functions for the $h$-polynomials of \sti.
\end{abstract}

\maketitle

{\openup \myspread  

\section{Introduction}\label{sec:intro}
\subsection{Overview}
Let $S_n$ denote the symmetric group of permutations of
$[n]:=\{1,\ldots,n\}$. For each $w\in S_n$, a \emph{descent}
is an index $i$ ($1\le i<n$) such that $w_i>w_{i+1}$. We let
\[
d(w):=|\{ i \in [n-1]: w_i > w_{i+1} \}|
\]
denote the number of descents in $w$.
The corresponding generating function
\begin{equation}\label{eq:classic}
A_{n-1}(t) := \sum_{w\in S_n}t^{\des(w)}
\end{equation}
is known as an \emph{Eulerian polynomial}, although this definition
differs from the classical one by a power of $t$. Some interesting
features of the Eulerian polynomials include the facts that they have
symmetric and unimodal coefficients and are known to have all real
roots.

More generally, if $W$ is any finite Coxeter group with simple
reflections $s_1,\dots,s_n$ (such as the symmetric group $S_{n+1}$
with simple transpositions $s_i=(i,i+1)$), then a descent in some
$w\in W$ may be defined as an index $i$ such that $\ell(ws_i)<\ell(w)$,
where $\ell(w)$ denotes the minimum length of an expression for $w$
as a product of simple reflections. Thus there is an analogous
\emph{$W$-Eulerian polynomial}
\[
W(t):=\sum_{w\in W}t^{d(w)},
\]
where $d(w)$ is defined to be the number of descents in $w$. Note that
as a Coxeter group, $S_n$ is often denoted $A_{n-1}$, so this notation
is consistent with~\eqref{eq:classic}.

Like the classical Eulerian polynomials, the $W$-Eulerian polynomials
are known to have symmetric and unimodal coefficients. An elegant
explanation of this fact may be based on a topological interpretation
of $W(t)$ as the $h$-polynomial of the Coxeter complex of $W$.
Since every (finite) Coxeter complex is realizable as the boundary
complex of a simplicial polytope, the symmetry and unimodality of the
coefficients of $W(t)$ may thus be seen as a consequence of the
$g$-theorem (e.g., see Section III.1 of~\cite{Stan}).

Recently, several authors
(see for example \cite{Branden, Gal, PostReinWill, Stembridge})
have identified interesting classes of simplicial complexes
whose $h$-polynomials have expansions of the form
\[
h(t) = \sum_{0\le i \le n/2} \gamma_i t^i(1+t)^{n-2i},
\]
where the coefficients $\gamma_i$ are nonnegative. It is easy to
see that each summand in this expansion has symmetric and unimodal
coefficients centered at $n/2$, and thus any $h$-polynomial
with a nonnegative ``$\gamma$-vector'' in this sense necessarily
has symmetric and unimodal coefficients. In these terms,
the $h$-polynomials of all finite Coxeter complexes
(i.e., the $W$-Eulerian polynomials) are known
to have nonnegative $\gamma$-vectors~\cite{Stembridge}.

Another feature of $\gamma$-nonnegativity is that it is a necessary
condition for a polynomial to have all real roots,
given that the polynomial has nonnegative symmetric coefficients.
In this direction, Brenti~\cite{Brenti} has conjectured that
the $W$-Eulerian polynomials have all real roots, a result that
remains unproved only for the groups $W=D_n$.

In this paper, we study a family of Eulerian-like polynomials
associated to irreducible affine Weyl groups. These ``affine''
Eulerian polynomials may be defined as generating functions for
``affine descents'' over the corresponding finite Weyl group.
An affine descent is similar to an ordinary descent in a Weyl
group, except that the reflection corresponding to the highest root
may also contribute a descent, depending on its effect on length.

The affine Eulerian polynomials have a number of interesting properties
similar to those of the ordinary $W$-Eulerian polynomials. In particular,
we show that they have nonnegative $\gamma$-vectors
(Theorem~\ref{thm:gamnonn}), and conjecture that all of their
roots are real. Perhaps the most interesting similarity is that each
affine Eulerian polynomial is the $h$-polynomial of a naturally
associated relative cell complex (Theorem~\ref{thm:torus}).

To describe this complex, one should start with an irreducible affine
Coxeter arrangement. Such an arrangement induces a simplicial
decomposition of the ambient space; by taking the quotient of this
space by the translation subgroup of the associated affine Weyl group,
one obtains a torus decomposed into simplicial cells.
We refer to this cell complex as the \emph{\st} in recognition
of the work of Steinberg, who gave a beautiful proof of Bott's formula
for the Poincar\'e series of an affine Weyl group by analyzing the action
of the finite Weyl group on the homology of this complex in two different
ways (see Section~3 of~\cite{Steinberg}). In fact, Steinberg also allows
the possibility of twisting the entire construction by an automorphism,
but we will not consider this variation here.

It is important to note that the \st{} is not a simplicial complex
(distinct cells may share the same set of vertices), but it is at
least a Boolean cell complex in the sense that all lower intervals
in the partial ordering of cells are Boolean algebras.\footnote{We
thank V.~Welker for bringing this to our attention.}
For further information about Boolean complexes, see~\cite{Stan2}
and the references cited there.

For our purposes, it is essential to omit the empty cell of dimension
$-1$ from the \st; we refer to the resulting relative complex as
the \emph{\rst}. It is this complex whose $h$-polynomial is the
corresponding affine Eulerian polynomial; i.e., the generating function
for affine descents.

It is noteworthy that affine descents in finite Weyl groups
were first introduced by Cellini~\cite{Cellini} in a construction of
a variant of Solomon's descent algebra, and developed further for
the groups of type $A$ and $C$ in several follow-up papers on
``cyclic descents'' by Cellini~\cite{Cellini2,Cellini3},
Fulman~\cite{Fulman0,Fulman}, and Petersen~\cite{Petersen}.
In very recent work, Lam and Postnikov~\cite{LamPost} study a
weighted count of affine descents (the ``circular descent number'')
that coincides with an ordinary count (only) in type $A$.

\subsection{Organization}
The paper is structured as follows. Section~\ref{sec:prelim}
introduces the necessary definitions, including details of the
construction of the \st. In Section~\ref{sec:affE} we show that the
affine Eulerian polynomials are the $h$-polynomials of \rsti{}
(Theorem~\ref{thm:torus}). Although we do not know of any conceptual
topological explanation for the nonnegativity of the $h$-vector,
we do show that \rsti{} are partitionable
(Remark~\ref{rmk:partition}); this is a weak analogue of shellability
that implies $h$-nonnegativity.

In Section~\ref{sec:rrgp}, we present our second main result;
namely, that the affine Eulerian polynomials have nonnegative
$\gamma$-vectors (Theorem~\ref{thm:gamnonn}). As a corollary,
it follows that the $h$-vectors of \rsti{} are symmetric and unimodal.
In this section, we also present evidence supporting our conjecture
that all roots of affine Eulerian polynomials are real. The proof
of Theorem~\ref{thm:gamnonn} is case-by-case, and relies on
combinatorial expansions for the $\gamma$-vectors of affine Eulerian
polynomials for the classical Weyl groups that we provide in
Section~\ref{sec:comb}. In this latter section, we also provide
combinatorial expansions for the flag $h$-polynomials of \rsti,
one of which suggests the possibility that a natural class of
(reduced) polyhedral tori may have nonnegative $cd$-indices
(see Question~\ref{quest:cd}).

In Section~\ref{sec:ident}, we present three unexpected identities
relating ordinary and affine Eulerian polynomials (two new, one old),
and use these to derive exponential generating functions for the affine
Eulerian polynomials for each classical series of Weyl groups.

\newpage 
\section{Preliminaries}\label{sec:prelim}
\subsection{Finite and affine Weyl groups}

We assume the reader is familiar with the basic theory of reflection
groups. We follow the notational conventions of~\cite{Humphreys}.

Let $\Phi$ be a \xtal{} root system embedded in a real
Euclidean space $V$ with inner product $\br{\cdot\,{,}\,\cdot}$.
For any root $\beta \in \Phi$,
let $H_\beta := \{\lambda \in V:\br{\lambda, \beta} = 0\}$
be the hyperplane orthogonal to $\beta$ and let $s_\beta$ denote
the orthogonal reflection through $H_\beta$. Fix a set of simple
roots $\Delta=\{\alpha_1,\dots,\alpha_n\} \subset \Phi$,
and let $S=\{s_1,\dots,s_n\}$ denote the corresponding set
of simple reflections. The latter generates a finite Coxeter
group $W$ (a Weyl group).

Unless stated otherwise, we always assume that $\Phi$ and $W$
are irreducible.

For convenience, we assume that $\Delta$ spans $V$.

Having fixed a choice of simple roots, every root $\beta$ either
belongs to the nonnegative span of the simple roots and is designated
\emph{positive}, or else belongs to the nonpositive span of the
simple roots and is designated \emph{negative}.
We write $\beta>0$ or $\beta<0$ accordingly.

The affine Weyl group $\aff{W}$ is generated by reflections
$s_{\beta,k}$ through the affine hyperplanes
\[
H_{\beta,k} := \{\lambda \in V:\br{\lambda, \beta} = k\}
\qquad(\beta\in\Phi,\ k\in\Z).
\]
Alternatively, one may construct $\aff{W}$ as the semidirect product
$W\ltimes\Z\Phi^\vee$, where $\Z\Phi^\vee$ denotes the lattice
generated by all co-roots $\beta^\vee=2\beta/\br{\beta,\beta}$
($\beta\in\Phi$), acting on $V$ via translations.

Given that $\Phi$ is irreducible, it has a unique highest root $\hr$,
and it is well-known that $\aff{W}$ is generated by
$\aff{S}:=S\cup\{s_{\hr,1}\}$ and that $(\aff{W},\aff{S})$ is an
irreducible Coxeter system.

Note that $\aff{W}$ depends on the underlying root system $\Phi$
(not merely $W$), so we are committing an abuse of notation.
For example, $B_n$ and $C_n$ are isomorphic as Coxeter systems,
but the affine groups $\aff{B}_n$ and $\aff{C}_n$ are not isomorphic
as Coxeter systems for $n\ge3$.

\subsection{Coxeter complexes}
The hyperplanes $H_\beta$ ($\beta\in\Phi$) induce a partition of $V$
into a complete $W$-symmetric fan of simplicial cones.
By intersecting this fan with the unit sphere in~$V$,
one obtains a topological realization of the Coxeter
complex $\Sigma(W)$.
The action of $W$ on chambers (maximal cones) in the fan is
simply transitive, and the choice of simple roots $\Delta$ is
equivalent to designating a dominant chamber; namely,
\[
C_\mtset:=\{\lambda\in V:\br{\lambda,\alpha}>0
  \text{ for all }\alpha\in\Delta\}.
\]
The closure of the dominant chamber is a fundamental domain for
the action of $W$ on $V$, and thus every cone in the fan has the
form $wC_J$ ($w\in W$, $J\subseteq[n]$), where
\[
C_J:=\{\lambda\in V: \br{\lambda, \alpha_j} = 0 \text{ for } j\in J,
  \ \br{\lambda, \alpha_j} > 0 \text{ for } j\in [n]\setminus J\}.
\]
Notice that the rays (1-dimensional cones) have the form $wC_J$
where $J=[n]\setminus\{j\}$ for some~$j$. If we assign color $j$
to all such rays, we obtain a \emph{balanced} coloring of $\Sigma(W)$;
i.e., every maximal face (chamber) has exactly one vertex (extreme ray)
of each color.

Similarly, the affine hyperplanes $H_{\beta,k}$ ($\beta\in\Phi$, $k\in\Z$)
may be used to partition $V$ into a $\aff{W}$-symmetric simplicial
complex that is isomorphic to the Coxeter complex $\Sigma(\aff{W})$.
By abuse of notation, we will identify $\Sigma(\aff{W})$ with this
particular geometric realization. The action of $\aff{W}$ on alcoves
(maximal simplices) is simply transitive, and the fundamental alcove
\[
A_\mtset:= C_\mtset \cap \{ \lambda \in V: \br{\lambda, \hr} < 1\}
\]
is tied to the choice of $\aff{S}$ in the sense that the
$\aff{W}$-stabilizer of every point in the closure of $A_\mtset$
(a fundamental domain) is generated by a proper subset of $\aff{S}$.
We index the faces of $A_\mtset$ by subsets of
$[0,n] := \{ 0,1,\ldots,n\}$ so that the $J$-th face is
\[
A_J:=\begin{cases}
  \hbox to 31pt{\hfill $C_J$} \cap \{ \lambda \in V: \br{\lambda, \hr} < 1\}
    &\text{if $0\notin J$},\\
  C_{J\setminus\{0\}} \cap \{ \lambda \in V: \br{\lambda, \hr} = 1\}
    &\text{if $0\in J$.}
\end{cases}
\]
Note that $A_J$ is the empty face when $J=[0,n]$.

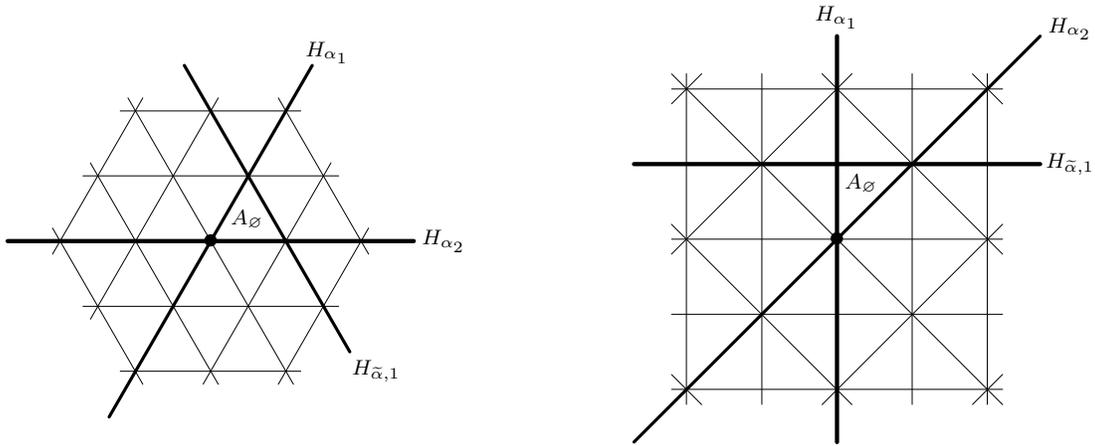
\begin{figure}
\begin{minipage}[c]{0.5\linewidth}
\vspace{11pt}
\hspace{-0.5\linewidth} 
\begin{xy}
0;<1cm,0cm>:<.5cm,\halfrootthree cm>::
(-4.5,0)*{},
(-.2,5); (2.2,5) **@{-},
(-.2,4); (3.2,4) **@{-},
(-.7,3); (4.7,3) **@{-},
(-.7,3); (4.7,3) **[|(4)]@{-},
(.8,2); (4.2,2) **@{-},
(1.8,1); (4.2,1) **@{-},
(0,2.8); (0,5.2) **@{-},
(1,1.8); (1,5.2) **@{-},
(2,.3); (2,5.7) **@{-},
(2,.3); (2,5.7) **[|(3)]@{-},
(3,.8); (3,4.2) **@{-},
(4,.8); (4,3.2) **@{-},
(-.2,3.2); (2.2,.8) **@{-},
(-.2,4.2); (3.2,.8) **@{-},
(-.2,5.2); (4.2,.8) **@{-},
(.3,5.7); (4.7,1.3) **@{-},
(.3,5.7); (4.7,1.3) **[|(3)]@{-},
(1.8,5.2); (4.2,2.8) **@{-},
(2.1,5.9)*{\scriptstyle H_{\alpha_1}},
(5.1,3)*{\scriptstyle H_{\alpha_2}},
(5.2,1)*{\scriptstyle H_{\hr,1}},
(2.33,3.33)*{\scriptstyle A_\mtset},
(2,3)*{\bullet}
\end{xy}
\end{minipage}%
\begin{minipage}[c]{0.5\linewidth}
\hspace{-0.5\linewidth} 
\begin{xy}
0;<1cm,0cm>:
(-6,0)*{},
(-.2,0); (4.2,0) **@{-},
(-.2,1); (4.2,1) **@{-},
(-.2,2); (4.2,2) **@{-},
(-.7,3); (4.7,3) **@{-},
(-.7,3); (4.7,3) **[|(4)]@{-},
(-.2,4); (4.2,4) **@{-},
(0,-.2); (0,4.2) **@{-},
(1,-.2); (1,4.2) **@{-},
(2,-.7); (2,4.7) **@{-},
(2,-.7); (2,4.7) **[|(4)]@{-},
(3,-.2); (3,4.2) **@{-},
(4,-.2); (4,4.2) **@{-},
(-.2,3.8); (.2,4.2) **@{-},
(-.2,1.8); (2.2,4.2) **@{-},
(-.7,-.7); (4.7,4.7) **@{-},
(-.7,-.7); (4.7,4.7) **[|(3)]@{-},
(1.8,-.2); (4.2,2.2) **@{-},
(3.8,-.2); (4.2,.2) **@{-},
(-.2,.2);  (.2,-.2) **@{-},
(-.2,2.2); (2.2,-.2) **@{-},
(-.2,4.2);  (4.2,-.2) **@{-},
(1.8,4.2); (4.2,1.8) **@{-},
(3.8,4.2); (4.2,3.8) **@{-},
(2,4.95)*{\scriptstyle H_{\alpha_1}},
(5.1,3)*{\scriptstyle H_{\hr,1}},
(5.1,4.8)*{\scriptstyle H_{\alpha_2}},
(2.33,2.73)*{\scriptstyle A_{\mtset}},
(2,2)*{\bullet}
\end{xy}
\end{minipage}
\caption{Portions of the Coxeter complexes
  for $\aff{A}_2$ and $\aff{C}_2$.}
\label{fig:affA2C2}
\end{figure}

The Coxeter complexes for $\aff{A}_2$ and $\aff{C}_2$
are illustrated in Figure~\ref{fig:affA2C2}.

Since the closure of $A_\mtset$ is a fundamental domain for the action
of $\aff{W}$, each cell in this complex has the form
$\mu+wA_J$ $(\mu\in\Z\Phi^\vee$, $w\in W$, $J \subseteq [0,n])$.
In particular, the vertices of $\Sigma(\aff{W})$ are of the form
$\mu+wA_{\{j\}^c}$, where $J^c:= [0,n]\setminus J$.
If we assign color $j$ to each of the vertices $\mu+wA_{\{j\}^c}$,
then the vertices of the cell $\mu+wA_J$ are assigned color-set $J^c$
(without repetitions), so this coloring is balanced.

\begin{rmk}\label{rmk:stW}
If $\Phi$ and $W$ are reducible, then the affine hyperplanes
$H_{\beta,k}$ may still be used to partition $V$ into a cell complex,
but the result is not a geometric realization of the Coxeter complex
of~$\aff{W}$. Indeed, the cells of this complex are products of
simplices, whereas the Coxeter complex of every Coxeter system is
simplicial.
\end{rmk}

\subsection{Flag $f$-vectors and $h$-vectors}\label{sec:fhpoly}
Let $\Sigma$ be a finite set of simplices (or abstractly,
a hypergraph) that is properly colored; i.e., the vertices of $\Sigma$
have been assigned colors from some index set, say $[0,n]$, so that
no simplex has two vertices with the same color. The main examples
we have in mind are balanced simplicial (or more generally, Boolean)
complexes.

A basic combinatorial invariant of $\Sigma$ that carries significant
algebraic and topological information (e.g., see the discussion
in Section~III.4 of~\cite{Stan}) is the \emph{flag $h$-vector}.
The components of the flag $h$-vector are the quantities
\begin{equation}\label{eq:hJdef}
h_J(\Sigma) := \sum_{I \subseteq J} (-1)^{|J\setminus I|}f_I(\Sigma)
  \qquad(J\subseteq[0,n]),
\end{equation}
where $f_I(\Sigma)$ denotes the number of simplices in $\Sigma$
whose vertices have color-set~$I$.

The quantities $f_J(\Sigma)$ for $J\subseteq[0,n]$ are collectively
referred to as the \emph{flag $f$-vector} of $\Sigma$.

The corresponding generating functions
\begin{align*}
f(\Sigma; t_0,\ldots,t_n) &:= \sum_{J\subseteq [0,n]}
  f_J(\Sigma)\prod_{j \in J} t_j,\\
h(\Sigma; t_0,\ldots,t_n) &:= \sum_{J\subseteq [0,n]}
  h_J(\Sigma)\prod_{j \in J} t_j
\end{align*}
are known as the \emph{flag $f$-polynomial} and \emph{flag $h$-polynomial}
of $\Sigma$. The more familiar ordinary $f$-polynomial and
$h$-polynomial may be obtained via the specializations
\begin{align*}
f(\Sigma; t)&:= f(\Sigma; t,\ldots,t)
  =\sum_{J\subseteq[0,n]}f_J(\Sigma)t^{|J|},\\
h(\Sigma; t)&:= h(\Sigma; t,\ldots,t)
  =\sum_{J\subseteq[0,n]}h_J(\Sigma)t^{|J|}.
\end{align*}
The coefficients of these polynomials yield
the (ordinary) $f$-vector and $h$-vector of $\Sigma$.

Note that \eqref{eq:hJdef} implies
\begin{equation}\label{eq:hfid}
h(\Sigma; t_0,\ldots,t_n) = (1-t_0) \cdots (1-t_n)
  f\Bigl(\Sigma; \frac{t_0}{1-t_0},\ldots, \frac{t_n}{1-t_n}\Bigr),
\end{equation}
and hence $h(\Sigma;t) = (1-t)^{n+1}f(\Sigma;t/(1-t))$.

\subsection{The \st}
As the translation subgroup of $\aff{W}$, the co-root lattice
$\Z\Phi^\vee$ acts as a group of color-preserving automorphisms of
the affine Coxeter complex~$\Sigma(\aff{W})$. Letting $T$ denote
the $n$-torus $V/\Z\Phi^\vee$, it follows that the image of
$\Sigma(\aff{W})$ under the natural map $V\to T$ is a balanced
Boolean complex, denoted $\stor{W}$. That is,
\[
\stor{W}=\Sigma(\aff{W})/\Z\Phi^\vee.
\]
As explained in the introduction, we refer to $\stor{W}$ as the \emph{\st}.
We also define the \emph{\rst}, denoted $\rstor{W}$, to be the relative
complex obtained by deleting the empty simplex of dimension $-1$
from $\stor{W}$.

Note that these are finite complexes;
there is one maximal cell $wA_\mtset+\Z\Phi^\vee$ for each $w\in W$.

There is an alternative way to construct the \st{} that starts with the
observation that the 0-colored vertices in $\Sigma(\aff{W})$ are the
members of $\Z\Phi^\vee$. Since every alcove $A$ has a unique 0-colored
vertex, one may translate $A$ via $\Z\Phi^\vee$ to a unique alcove that
has the origin as a vertex; i.e., to one of the alcoves in the $W$-orbit
of $A_\mtset$. The closure of this set of alcoves is the $W$-invariant
convex polytope
\[
P_\Phi=\{\lambda\in V:
  -1\le\br{\lambda,\beta}\le 1\text{ for all }\beta\in\Phi\},
\]
and the \st{} is obtained by identifying the maximal opposite
faces of $P_\Phi$.

\begin{figure}
\centering
\text{\vbox to1.21in{\vfill\hbox{\includegraphics{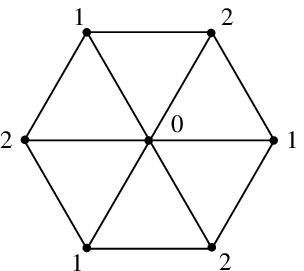}}\vfill}}
\hskip 0.15\linewidth
\includegraphics{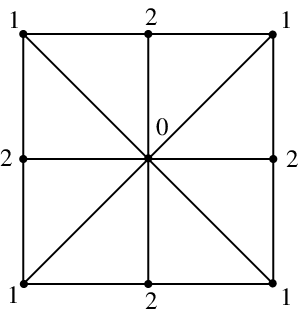}
\caption{The \sti{} for $\aff{A}_2$ and $\aff{C}_2$.}
\label{fig:tori}
\end{figure}

\begin{exm}\label{exm:A2}
The \st{} for $\aff{A}_2$ is a hexagon with opposite sides identified,
decomposed into six triangles, nine edges, and three vertices.
See Figure~\ref{fig:tori}. It has flag $f$-polynomial
\[
f(\storr{A}{2};t_0,t_1,t_2)
  =1+t_0+t_1+t_2+3t_0t_1+3t_0t_2+3t_1t_2+6t_0t_1t_2.
\]
Using~\eqref{eq:hfid} to compute the flag $h$-polynomial, we find
\[
h(\storr{A}{2};t_0,t_1,t_2)
  =1 + 2 t_0t_1 + 2 t_0t_2 + 2t_1t_2 - t_0 t_1 t_2.
\]
On the other hand, the \rst{} lacks the empty face, so its
flag $f$-polynomial omits the constant term and we find
\[
h(\rstorr{A}{2};t_0,t_1,t_2)
  =t_0 + t_1 + t_2 + t_0t_1 + t_0t_2 + t_1t_2.
\]
Specializing, we see that the \rst{} has ordinary $f$-polynomial
$3t+9t^2+6t^3$, and ordinary $h$-polynomial $3t+3t^2$.
\end{exm}

\begin{exm}\label{exm:C2}
The \st{} for $\aff{C}_2$ (or the isomorphic $\aff{B}_2$) is a square
with opposite sides identified, decomposed into eight triangles,
twelve edges, and four vertices as in Figure~\ref{fig:tori}.
The \rst{} has flag $f$-polynomial
\[
f(\rstorr{C}{2};t_0,t_1,t_2)
  =t_0+t_1+2t_2+4t_0t_1+4t_0t_2+4t_1t_2+8t_0t_1t_2,
\]
and (again via~\eqref{eq:hfid}) flag $h$-polynomial
\[
h(\rstorr{C}{2};t_0,t_1,t_2)
  =t_0+t_1+2t_2+2t_0t_1+t_0t_2+t_1t_2.
\]
As in the previous example, it is easy to check that the ordinary
and flag $h$-polynomials of the unreduced \st{} have (some) negative
coefficients.
\end{exm}

\subsection{Affine descents}
We define a root $\beta$ to be
negative with respect to $w\in W$ if $w\beta<0$.
The positive roots that are negative with respect to $w$
are known as \emph{inversions}.
If $\ell(w)$ denotes the minimum length of an expression for $w$
as a product of simple reflections, then $\beta$ is negative with
respect to $w$ if $\ell(ws_\beta)<\ell(w)$ (for $\beta>0$)
or $\ell(ws_\beta)>\ell(w)$ (for $\beta<0$).

A simple root that is negative with respect to $w$ is said to be a
(right) \emph{descent}, and the descent set of $w$, denoted $D(w)$,
records the corresponding set of indices. Thus,
\[
\Des(w) = \{ j \in [n] : w\alpha_j < 0 \}
  =\{j\in[n]:\ell(ws_j)<\ell(w)\}.
\]
We let $\des(w):= |\Des(w)|$ denote the number of descents in $w$.

As noted in the introduction, the $W$-Eulerian polynomial is the
$h$-polynomial of the Coxeter complex $\Sigma(W)$. That is,
\[
W(t)=\sum_{w\in W}t^{\des(w)}=h(\Sigma(W);t).
\]
More generally, the generating function for descent \emph{sets}; namely,
\[
W(t_1,\dots,t_n):=\sum_{w\in W}\prod_{j\in\Des(w)}t_j
\]
is the flag $h$-polynomial of $\Sigma(W)$
(e.g., see the discussion at the end of Section~III.4 in~\cite{Stan}).

Extending these concepts, set $\alpha_0:=-\hr$ (the lowest root),
and let $s_0=s_{\hr}$ denote the corresponding reflection in~$W$.
We define the \emph{affine descent} set of $w$, denoted $\aDes(w)$,
to be the set of indices of roots in $\Delta_0:=\Delta\cup\{\alpha_0\}$
that are negative with respect to $w$. Thus,
\[
\aDes(w) = \{ j \in [0,n] : w\alpha_j < 0 \}
  =\begin{cases}
  \Des(w)\cup\{0\}&\text{if $\ell(ws_0)>\ell(w)$},\\
  \cbox{59}{\Des(w)}&\text{if $\ell(ws_0)<\ell(w)$.}\
  \end{cases}
\]
We let $\ades(w):= |\aDes(w)|$ denote the number of affine descents in $w$.

Note that only the identity element of $W$ has an empty descent set
(but has an affine descent at~0), and only the longest element $w_0$
has a full descent set (i.e., $\Des(w_0)=[n]$) but does not have
an affine descent at~0.  Thus $1\le\ades(w)\le n$ for all $w\in W$.

\section{Affine Eulerian polynomials}\label{sec:affE}
We let $\aff{W}(t_0,\dots,t_n)$ and $\aff{W}(t)$ denote the
respective generating functions for affine descent sets and numbers
of affine descent sets; i.e.,
\begin{align}
\aff{W}(t_0,\dots,t_n)
  &:=\sum_{w\in W}\prod_{j\in\aDes(w)}t_j,\label{eq:afep}\\
\aff{W}(t):=\aff{W}(t,\dots,t)
  &\phantom{:}=\sum_{w\in W}t^{\ades(w)}.\label{eq:aoep}
\end{align}
We refer to these as \emph{multivariate} and
\emph{univariate affine Eulerian polynomials}.

\begin{thm}\label{thm:torus}
If $\aff{W}$ is an irreducible affine Weyl group,
then the flag $h$-polynomial of the corresponding \rst{}
is the multivariate $\aff{W}$-Eulerian polynomial; i.e.,
\begin{equation}\label{eq:heqW}
h(\rstor{W};t_0,\dots,t_n)=\aff{W}(t_0,\dots,t_n).
\end{equation}
In particular, for all $J\subseteq[0,n]$, we have
\begin{align}
f_J(\rstor{W}) &= |\{ w \in W: \aDes(w)\subseteq J\}|,\label{eq:flagf}\\
h_J(\rstor{W}) &= |\{ w \in W: \aDes(w)=J\}|.\label{eq:flagh}
\end{align}
Furthermore,
\begin{equation}\label{eq:affEu}
\aff{W}(t_0,\dots,t_n)=\sum_{J\subsetneq[0,n]}
  \frac{|W|}{|W_J|}\prod_{j\in J}(1-t_j)\prod_{j\notin J}t_j,
\end{equation}
where $W_J$ denotes the (not necessarily parabolic) subgroup
of $W$ generated by $\{s_j:j\in J\}$.
\end{thm}

Of course it follows immediately that the ordinary $h$-polynomial
of the \rst{} is the corresponding univariate affine
Eulerian polynomial; i.e.,
\[
h(\rstor{W}; t) = \aff{W}(t)=\sum_{w \in W} t^{\ades(w)}.
\]

\begin{cor}\label{cor:sym}
The flag $h$-vector of the \rst{} $\rstor{W}$ satisfies the
generalized Dehn-Sommerville equations; that is,
for all $J\subseteq[0,n]$, we have
\[
h_J(\rstor{W})=h_{J^c}(\rstor{W}).
\]
In particular, the $\aff{W}$-Eulerian polynomial is symmetric: $\aff 
{W}(t)=t^{n+1}\aff{W}(1/t)$.
\end{cor}

\begin{proof}
Recall that the longest element $w_0\in W$ is an involution that sends
all positive roots to negative roots. It follows that a root $\beta$
satisfies $w\beta<0$ if and only if $w_0w\beta>0$, and hence
\[
\aDes(w_0w) = [0,n] \setminus \aDes(w),
\]
for all $w\in W$. Now apply~\eqref{eq:flagh}.
\end{proof}

\begin{rmk}
The unreduced \st{} $\stor{W}$ has nearly the same flag $f$-vector
as its reduced counterpart, the only difference being
$f_\mtset(\stor{W})=1$ in place of $f_\mtset(\rstor{W}){=0}$.
However, as we noted in Example~\ref{exm:A2}, the $h$-polynomial need
not have symmetric or nonnegative coefficients in the unreduced case,
and is therefore of less interest.
\end{rmk}

The following lemma is the key to our proof of Theorem~\ref{thm:torus}.

\begin{lemma}\label{lem:coset}
If $\,\{\beta_i:i\in I\}$ is a set of simple roots for a
reflection subgroup $W'$ of $W$, then every coset in $W/W'$
has a unique member $w$ such that $w\beta_i>0$ for all $i\in I$.
\end{lemma}

\begin{proof}
Fix a dominant point $\lambda\in C_\mtset$ (i.e., $\br{\lambda,\alpha}>0$
for all roots $\alpha\in\Delta$), so that the $W$-orbit of $\lambda$ is
generic and the map $w\mapsto w^{-1}\lambda$ is a bijection between $W$
and the orbit $W\lambda$. Since
\[
w\beta_i>0\Leftrightarrow\br{\lambda,w\beta_i}>0
  \Leftrightarrow\br{w^{-1}\lambda,\beta_i}>0,
\]
we see that $w$ satisfies $w\beta_i>0$ for all $i\in I$ if and
only if $w^{-1}\lambda$ is dominant with respect to the simple
roots of $W'$. However, every $W'$-orbit has a unique dominant  
member, and the image of the coset $wW'$ under the bijection
is the $W'$-orbit of $w^{-1}\lambda$, so the result follows.
\end{proof}

\begin{rmk}
In the above lemma, it is interesting to note that by choosing the
simple roots of $W'$ so that they are positive relative to $\Phi$,
one may deduce that every coset of every reflection subgroup of $W$
has a unique element of minimum length. This is a familiar fact for
parabolic subgroups, but the less familiar general case also follows
from work of Dyer (see Corollary~3.4 of~\cite{Dyer}).
\end{rmk}

\begin{proof}[Proof of Theorem~\ref{thm:torus}.]
For each nonempty $J\subseteq[0,n]$, the set of cells of the affine
Coxeter complex with color-set $J$ is the $\aff{W}$-orbit of $A_{J^c}$.
These are the cells of the form $\mu+wA_{J^c}$
(for $\mu\in\Z\Phi^\vee$, $w\in W$), so
\[
\{wA_{J^c}+\Z\Phi^\vee:w\in W\}
\]
is the set of cells of the \rst{} with color-set $J$. However,
the $\aff{W}$-stabilizer of $A_{J^c}$ (or indeed, any subset of
the closure of the fundamental alcove) is generated by the subset
of $\aff{S}$ that fixes $A_{J^c}$. The $W$-image of this subgroup
(i.e., the $W$-stabilizer of $A_{J^c}+\Z\Phi^\vee$) is $W_{J^c}$,
the reflection subgroup of $W$ generated
by $\{s_j:j\in [0,n]\setminus J\}$, and therefore
\begin{equation}\label{eq:niceJ}
f_J(\rstor{W})=|W|/|W_{J^c}|.
\end{equation}
On the other hand, we have
\[
\{w\in W:\aDes(w)\subseteq J\}
  =\{w\in W:w\alpha_j>0\text{ for }j\in[0,n]\setminus J\}
\]
and every proper subset of $\Delta_0$ is the set of simple roots of
some root subsystem of $\Phi$ (this amounts to the fact that every
proper subset of the extended Dynkin diagram, which records the geometry
of $\Delta_0$, is the Dynkin diagram of a finite root system),
so Lemma~\ref{lem:coset} implies that $\{w\in W:\aDes(w)\subseteq J\}$
is a set of coset representatives for $W/W_{J^c}$. Hence,
\[
f_J(\rstor{W})=|W|/|W_{J^c}|=|\{w\in W:\aDes(w)\subseteq J\}|.
\]
Noting that $f_J(\rstor{W})=|\{w\in W:\aDes(w)\subseteq J\}|=0$
when $J=\mtset$, we obtain~\eqref{eq:flagf}.

To complete the proof, note that~\eqref{eq:flagf} implies
\eqref{eq:flagh} and hence~\eqref{eq:heqW} via inclusion-exclusion.
The latter allows one to deduce~\eqref{eq:affEu}
as a corollary of~\eqref{eq:hfid} and~\eqref{eq:niceJ}.
\end{proof}

\begin{table}[t]
\newcommand{\myspace}{\hline\vphantom{\vbox to 14pt{}}}
\newcommand{\mmyspace}{\hline\hline\vphantom{\vbox to 14pt{}}}
\begin{centering}
\begin{tabular}{ c | c }
$W$ & $\aff{W}(t)$ \\
\myspace
$B_3$ & $10t + 28t^2 + 10t^3$\\
\myspace
$B_4$ & $24t + 168t^2 + 168t^3 + 24t^4$\\
\myspace
$B_5$ & $54t + 904t^2 + 1924t^3 + 904t^4 + 54t^5$ \\
\myspace
$B_6$ & $116t + 4452t^2 + 18472t^3 + 18472t^4 + 4452t^5 + 116t^6$\\
\myspace
$B_7$ & $242t+20612t^2+157294t^3+288824t^4+157294t^5+20612t^6+242t^7$\\
\mmyspace
$D_4$ & $16t + 80t^2 + 80t^3 + 16t^4$\\
\myspace
$D_5$ & $44t + 464t^2 + 904t^3 + 464t^4 + 44t^5$\\
\myspace
$D_6$ & $104t + 2568t^2 + 8848t^3 + 8848t^4 + 2568t^5 + 104t^6$\\
\myspace
$D_7$ & $228t+13192t^2+79580t^3+136560t^4+79580t^5+13192t^6+228t^7$\\
\mmyspace
$E_6$ & $351t+5427t^2+20142t^3+20142t^4+5427t^5+351t^6$ \\
\myspace
$E_7$ & $4064t + 115728t^2 + 710112t^3 + 1243232t^4
   + 710112t^5 + 115728t^6+4064t^7$ \\
\myspace
\raise 8pt\hbox{$E_8$} & \vbox{\vglue 4pt\hbox to 373pt{%
    $157200t+9253680t^2+87417360t^3+251536560t^4$\hfill}
  \medskip\hbox to 373pt{\hfill
    ${}+251536560t^5+87417360t^6+9253680t^7+157200t^8$}}\\
\myspace
$F_4$ & $72t + 504t^2 + 504t^3 + 72t^4$\\
\myspace
$G_2$ & $6t + 6t^2$ \\
\hline
\end{tabular}
\bigskip
\end{centering}
\caption{Some affine Eulerian polynomials.}
\label{table}
\end{table}

It is easy to compute the affine Eulerian polynomials for the groups
of low rank via~\eqref{eq:affEu}. Some examples, including all of the
exceptional groups, are listed in Table~\ref{table}.

\begin{rmk}\label{rmk:partition}
Given $J\subsetneq[0,n]$, it follows from Lemma~\ref{lem:coset}
that each coset in $W/W_J$ has a unique representative $w$ such
that $\aDes(w)\cap J=\mtset$. Thus each cell of the \rst{} has the
form $F(w,J)=wA_J+\Z\Phi^\vee$ for some unique pair $(w,J)$ with
$\aDes(w)\cap J=\mtset$. Moreover, the cells of the form $F(w,*)$
are precisely the cells in the closure of $F(w,\mtset)$ that have
on their boundary the unique cell with color-set $\aDes(w)$;
namely, $F(w,\aDes(w)^c)$.  Thus the \rst{} is ``partitionable''
in the sense defined in Section~III.2 of~\cite{Stan}.
\end{rmk}

\begin{rmk}\label{rmk:fake}
If $W$ is an irreducible but non-\xtal{} finite reflection group,
such as $H_3$ or $H_4$, then there is no corresponding affine Weyl
group, and hence no \st. The root system still has a unique dominant
root, so one could define fake affine Eulerian polynomials
$\fake{W}(t_0,\dots,t_n)$ and $\fake{W}(t)$ analogous
to~\eqref{eq:afep} and~\eqref{eq:aoep}, using the
anti-dominant root in the role of~$\alpha_0$.
For the groups $H_3$ and $H_4$, one obtains
\begin{gather*}
\fake{H_3}(t)=26t + 68t^2 + 26t^3,\\
\fake{H_4}(t)=960t + 6240t^2 + 6240t^3 + 960t^4.
\end{gather*}
Although these polynomials have symmetric and unimodal coefficients
(and real roots), we do not know if $\fake{W}(t)$ is the $h$-polynomial
of some naturally associated Boolean complex.
\end{rmk}

\begin{rmk}
More generally, given any subset of roots
$\Psi=\{\beta_i:i\in I\}\subset\Phi$,
one could define a generalized descent set for $w\in W$ by setting
\[
\Des_\Psi(w):=\{i\in I:w\beta_i<0\},
\]
whether or not $\Phi$ is \xtal. Examining the proof of
Theorem~\ref{thm:torus}, one can see that the generating function for
these generalized descent sets would satisfy a formula similar
to~\eqref{eq:affEu} if for every $J\subseteq I$,
either $\{\beta_j:j\in J\}$ is the set of simple roots of some
finite root system (see Lemma~\ref{lem:coset}),
or $\{w\in W:w\beta_j>0\text{ for all $j\in J$}\}$ is empty.
Applying this criterion to pairs $i,j\in I$, this forces the angle
between $\beta_i$ and $\beta_j$ to be $(1-1/m)\pi$ for some
integer $m\ge2$, or $\beta_i=-\beta_j$ (i.e., $m=\infty$).
Since the matrix $\br{\beta_i,\beta_j}$ is necessarily positive
semidefinite, it follows from the theory of reflection groups
that (up to normalization) $\Psi$ must be the simple roots of some
root subsystem, or is an extension of the simple roots by the lowest
root of some \xtal{} root subsystem, or is an orthogonal
disjoint union of such sets (e.g., see Section~2.7 of~\cite{Humphreys}).
In particular, the identity in~\eqref{eq:affEu} is not valid for
the fake affine Eulerian polynomials discussed in the previous remark.
\end{rmk}

\section{Real roots, $\gamma$-vectors, and unimodality}\label{sec:rrgp}
The following is a companion to Brenti's conjecture~\cite{Brenti} that
the roots of all (ordinary) Eulerian polynomials $W(t)$ are real.

\begin{conj}\label{conj:rr}
The roots of all affine Eulerian polynomials $\aff{W}(t)$ are real.
\end{conj}

To complete a proof of this conjecture, we claim that it suffices to
consider only the groups $\aff{B}_n$ and $\aff{D}_n$.
Indeed, it follows from observations of Fulman~\cite{Fulman0,Fulman}
and Petersen~\cite{Petersen} that $\aff{A}_n(t)$ and $\aff{C}_n(t)$
are both multiples of $A_{n-1}(t)$ (see also the discussion in
Section~\ref{sec:comb} below). Thus the conjecture for $\aff{A}_n$
and $\aff{C}_n$ follows from the fact that all roots of the
classical Eulerian polynomials are known to be real~\cite{Harper}.
Furthermore, using the data in Table~\ref{table}, it is easy to check
that the conjecture holds for the exceptional groups.

To collect supporting evidence for
the remaining groups $\aff{B}_n$ and $\aff{D}_n$,
we have determined explicit exponential generating functions
for the corresponding affine Eulerian polynomials
(see~Proposition~\ref{prp:aBDgf} below),
and used these to verify the conjecture for $n\le 100$.
In a similar way, we have also confirmed that all roots of $D_n(t)$ are
real (the only remaining open case of Brenti's conjecture) for $n\le 100$.

A further supporting result involves $\gamma$-vectors
in the sense of~\Branden~\cite{Branden} and Gal~\cite{Gal}.
To~explain, consider a polynomial satisfying $h(t)=t^mh(1/t)$.
It is clear that such a polynomial has a unique expansion
of the form
\[
h(t) = \sum_{0\le i \le m/2} \gamma_i t^i (1+t)^{m-2i}.
\]
We call $(\gamma_0,\gamma_1,\ldots)$
the \emph{$\gamma$-vector} of $h(t)$.

It is elementary to show that if $h(t)$ has symmetric, nonnegative
coefficients and all real roots, then it has a nonnegative
$\gamma$-vector (see Lemma~4.1 of~\cite{Branden}
or Section~1.4 of~\cite{Stembridge}).

Recall that $\aff{W}(t)$ is symmetric (Corollary~\ref{cor:sym}),
so it has a $\gamma$-vector.

\begin{thm}\label{thm:gamnonn}
The affine Eulerian polynomials $\aff{W}(t)$
have nonnegative $\gamma$-vectors.
\end{thm}

\begin{proof}
Given that we know Conjecture~\ref{conj:rr} holds for
$\aff{A}_n$, $\aff{C}_n$, and the exceptional affine Weyl groups,
it suffices to prove this result for $\aff{B}_n$ and $\aff{D}_n$.
In these cases, we have explicit combinatorial expansions
for $\aff{B}_n(t)$ and $\aff{D}_n(t)$ in Corollaries~\ref{cor:gamB}
and~\ref{cor:gamD} below that transparently imply the nonnegativity
of their $\gamma$-vectors.
\end{proof}

It would be interesting to have a conceptual (case-free) proof
of this result.

Any polynomial with a nonnegative $\gamma$-vector
has unimodal coefficients. Hence,

\begin{cor}
The affine Eulerian polynomials have unimodal coefficients.
\end{cor}

We remark that the $\gamma$-vectors of the Eulerian polynomials $W(t)$
are also known to be nonnegative, but the only existing proofs to date
are case-by-case~\cite{Chow,Stembridge}.

\section{Combinatorial expansions
  and $\gamma$-nonnegativity}\label{sec:comb}

In this section, we provide combinatorial expansions for the
affine Eulerian polynomials (both multivariate and univariate)
for the four infinite families of irreducible Weyl groups.
As corollaries, we will deduce the nonnegativity of the
$\gamma$-vectors for these polynomials.

\subsection{Type $A$}

Recall that the Weyl group $A_{n-1}$ may be identified with $S_n$,
the symmetric group of permutations of $[n]$, and the corresponding
root system is
\[
\{\eps_i-\eps_j:1\le i\ne j\le n\},
\]
where $\eps_1,\dots,\eps_n$ is the standard orthonormal basis of $\R^n$.
For the simple roots, we choose $\alpha_i = \eps_{i+1}-\eps_i$
($1\le i<n$). With respect to this choice, the simple reflection $s_i$
transposes $i$ and $i+1$ (as a permutation) and interchanges
$\eps_i$ and $\eps_{i+1}$ (as a reflection acting on $\R^n$).

The positive roots are $\eps_i-\eps_j$ for $i>j$,
and the lowest root $\alpha_0=-\hr$ is $\eps_1-\eps_n$.

We write permutations in one-line form $w=w_1w_2\cdots w_n$,
where $w_i=w(i)$. In these terms, a root $\eps_i-\eps_j$ is
negative with respect to a permutation $w$ if and only if $w_j>w_i$.
In particular, $\Des(w)=\{i\in[n-1]:w_i>w_{i+1}\}$ is the usual
descent set of a permutation. Also, an extra ``affine'' descent
occurs at $0$ if and only if $w_n>w_1$,
so $\aDes(w)=\{i\in[0,n-1]:w_i>w_{i+1}\}$,
using the convention $w_0=w_n$.

For example, $\Des(25413)=\{2,3\}$ and $\aDes(25413)=\{0,2,3\}$.

\begin{prp}\label{prp:affA}
For $n\ge2$, we have
\[
\aff{A}_{n-1}(t_0,\dots,t_{n-1})=\sum_{j=0}^{n-1}
  t_j A_{n-2}(t_{j+1},\ldots,t_{n-1},t_0,\ldots,t_{j-2}).
\]
\end{prp}

\begin{proof}
Let $c=23\cdots n1$ (an $n$-cycle in $A_{n-1}$), and
note that one may obtain the affine descent set
of $wc=w_2\cdots w_nw_1$ by a cyclic shift of the affine descent
set of $w\in A_{n-1}$; i.e.,
\[
\aDes(wc)=\{ i-1 : i\in\aDes(w)\}\mod n.
\]
Each coset of the cyclic subgroup $\br{c}$ has a unique representative
$w$ such that $w_n=n$, and this set of representatives is in bijection
with $A_{n-2}$. For each coset representative $w$,
we have $0\in\aDes(w)$, and the remaining affine descents coincide
with the ordinary descents of the corresponding member of $A_{n-2}$.
Thus, the generating function for the affine descent sets of these
coset representatives is $t_0A_{n-2}(t_1,\dots,t_{n-2})$,
and the generating function corresponding to elements of the
form $wc^{-j}$ is obtained by substituting $t_i\to t_{i+j}$
(subscripts modulo~$n$).
\end{proof}

It follows that the univariate affine Eulerian
polynomials of type $A$ are multiples of classical Eulerian
polynomials, as noted previously by Fulman~\cite{Fulman0}
and Petersen~\cite{Petersen}.

\begin{cor}\label{cor:cycA}
For $n\ge1$, we have $\aff{A}_n(t)=(n+1)tA_{n-1}(t)$.
\end{cor}

\subsection{Type $C$}\label{sec:Cgam}
The root system of the Weyl group $C_n$ has the form
\[
\{\pm2\eps_i:1\le i\le n\}\cup
  \{\pm\eps_i\pm\eps_j:1\le j<i\le n\},
\]
and $C_n$ acts as a group of permutations
of $\{\pm\eps_1,\ldots,\pm\eps_n\}$. More explicitly,
if we identify $\pm i$ with $\pm\eps_i$, then $C_n$ may be viewed
as the group of permutations of $\pm[n]=\{\pm1,\dots,\pm n\}$
such that $w(-i)=-w(i)$ for all $i$.
The simple roots may be chosen so that $\alpha_1 = 2\eps_1$
and $\alpha_i = \eps_{i} - \eps_{i-1}$ for $2\le i\le n$.
With respect to this choice, the positive roots are $2\eps_i$
for all $i$ and $\eps_i\pm\eps_j$ for all $i>j$,
and the lowest root $\alpha_0=-\hr$ is $-2\eps_n$.

We write permutations $w\in C_n$ in one-line form $w=w_1\cdots w_n$,
where $w_i=w(i)$. In these terms, one can check that roots of the
form $\eps_i-\eps_j$ with $i>j$ are negative with respect to $w$ if
and only if $w_j>w_i$, whereas roots of the form $2\eps_i$ are
negative with respect to $w$ if and only if $w_i<0$. In particular,
the ordinary descent set is $\Des(w)=\{i\in[n]:w_{i-1}>w_i\}$,
using the convention $w_0=0$, and 0 is in the affine descent set
$\aDes(w)$ when $w_n>0$.

For example, if $w=23\bar{5}\bar{1}4$ (bars indicate
negative values), then $\aDes(w)=\{0,3\}$.

In the following formula for the multivariate $\aff{C}_n$-Eulerian
polynomial, it is more convenient to use $n+1$ in place
of $0$ to mark the extra affine descent, or equivalently,
set $t_0=t_{n+1}$. Note that by specializing this extra variable
(i.e., setting $t_0=t_{n+1}=1$), we recover Stembridge's formula for
the flag $h$-polynomial of the Coxeter complex $\Sigma(C_n)$
(Proposition~A.1 in~\cite{Stembridge}).

Below, we use $\chi(\cdot)$ as an indicator function:
$\chi(S)=1$ if $S$ is true; 0 if $S$ is false.

\begin{prp}\label{prp:flagC}
For $n\ge1$, we have
\[
\aff{C}_n(t_{n+1},t_1,\dots,t_n)=\sum_{u\in S_n}\prod_{i=1}^n
  \left(t_i^{\chi(u_{i-1}<u_i)}+t_{i+1}^{\chi(u_i>u_{i+1})}\right),
\]
using the convention $u_0=u_{n+1}=0$.
\end{prp}

\begin{proof}
Following the proof of Proposition~A.1 in~\cite{Stembridge},
each member of $C_n$ has the form $w=\sigma u$,
where $u\in S_n$ and $\sigma=(\sigma_1,\dots,\sigma_n)\in\Z_2^n$
(meaning that $w_i=\sigma_iu_i$). Given that $u_0=u_{n+1}=0$ and
that $n+1$ replaces 0 in $\aDes(w)$, we see that for $i=1,\dots,n+1$,
\begin{itemize}
\item if $u_{i-1}<u_i$, then $i\in\aDes(w)\Leftrightarrow\sigma_i=-1$,
\item if $u_{i-1}>u_i$, then $i\in\aDes(w)\Leftrightarrow\sigma_{i-1}=+1$.
\end{itemize}
Thus for each $i\in[n+1]$, there is a unique $j$ (depending on $u$)
such that the presence or absence of $i$ in $\aDes(w)$ is controlled by
the value of $\sigma_j$. More specifically, $\sigma_j$ controls the
presence of $j$ (if $u_{j-1}<u_j$) and $j+1$ (if $u_j>u_{j+1}$),
and nothing else. Furthermore, if we define
\begin{equation}\label{eq:csum}
c_j(u):=t_j^{\chi(u_{j-1}<u_j)}+t_{j+1}^{\chi(u_j>u_{j+1})},
\end{equation}
then $c_j(u)$ records the sum of the weights of the effects of
$\sigma_j=-1$ and $\sigma_j=+1$ on the affine descent set
of $\sigma u$. Since the effects of $\sigma_1,\dots,\sigma_n$ are
mutually independent, we conclude that
\[
\sum_{\sigma\in\Z_2^n\vphantom{\aDes}}\ \prod_{i\in\aDes(\sigma u)}
  \!\!t_i\ = \ c_1(u)\cdots c_n(u),
\]
and the result follows by summing over $u\in S_n$.
\end{proof}

\begin{rmk}
It is well-known that the Coxeter complex $\Sigma(C_n)$ is isomorphic
to the barycentric subdivision of an $n$-dimensional cube, and as
explained in Remark~A.3 of~\cite{Stembridge}, one may recognize
Stembridge's formula for the flag $h$-polynomial of $\Sigma(C_n)$ as
a disguised formula for the $cd$-index of the $n$-cube.
Similarly, the \st{} $\storr{C}{n}$ may be constructed from the
barycentric subdivision of an $n$-cube by identifying the opposite
maximal faces of the cube (recall Example~\ref{exm:C2}), and the
above formula may be reinterpreted as a nonnegative $cd$-index for
the reduced complex.
\end{rmk}

The above remark suggests the possibility of a more general result.
Given a tiling of $\R^n$ by lattice translates of a convex polytope $P$,
the quotient of $\R^n$ by the lattice may be viewed as an $n$-torus
decomposed into polyhedral cells. Lattice translates of the cells in
the barycentric subdivision of $P$ may be identified, thereby yielding
a ``polytorus'' with a well-defined flag $h$-vector, and our
findings here suggest that one should study the ``reduced polytorus''
obtained by deleting the empty face.

\begin{quest}\label{quest:cd}
Does every reduced polytorus have a nonnegative $cd$-index?
\end{quest}

It will be convenient for what follows to introduce three conventions
for counting peaks in a permutation $u\in S_n$; namely,
\begin{align*}
\pk(u)&:=|\{i\in[2,n-1]:u_{i-1}<u_i>u_{i+1}\}|,\\
\lpe(u)&:=|\{i\in[1,n-1]:u_{i-1}<u_i>u_{i+1}\}|,\\
\xpe(u)&:=|\{i\in[1,n]:u_{i-1}<u_i>u_{i+1}\}|,
\end{align*}
again using the convention $u_0=u_{n+1}=0$. We refer to these
quantities as the number of \emph{ordinary}, \emph{left},
and \emph{extended peaks} in $u$, respectively.

The following expansions show that $\aff{C}_n(t)$, $C_n(t)$,
and $A_{n-1}(t)$ have nonnegative $\gamma$-vectors. Part (b) is due
to Petersen (Proposition~4.15 in~\cite{Petersen}), and part (c) is
equivalent to an identity due to Foata and Sch\"utzenberger
(Th\'eor\`eme~5.6 of~\cite{FoSc}; see also Remark~4.8
of~\cite{StembridgeEnriched}).

\begin{cor}\label{cor:gamC}
For $n\ge1$, we have
\begin{itemize}\vspace{2pt}
\item[(a)] $\aff{C}_n(t)=(1/2)\sum_{u\in S_n}
  (4t)^{\xpe(u)}(1+t)^{n+1-2\xpe(u)}$,\vspace{6pt}
\item[(b)] $C_n(t)=\sum_{u\in S_n}
  (4t)^{\lpe(u)}(1+t)^{n-2\lpe(u)}$,\vspace{6pt}
\item[(c)] $A_{n-1}(t) = 2^{-(n-1)}\sum_{u \in S_n}
    (4t)^{\pk(u)}(1+t)^{n-1-2\pk(u)}$.
\end{itemize}
\end{cor}

\begin{proof}
(a) Proposition~\ref{prp:flagC} implies that
$\aff{C}_n(t_{n+1},t_1,\dots,t_n)=\sum_{u\in S_n}c_1(u)\cdots c_n(u)$,
where $c_i(u)$ is defined as in~\eqref{eq:csum}.
Specializing the variables so that $t_i\to t$ for all $i$,
one sees that
\begin{equation}\label{eq:crule}
c_i(u)\to\begin{cases}
  \cbox{20}{2t} &\text{if $u_{i-1}<u_i>u_{i+1}$ (a peak)},\\
  \cbox{20}{2}  &\text{if $u_{i-1}>u_i<u_{i+1}$ (a valley)},\\
  \cbox{20}{1+t}&\text{otherwise}.
\end{cases}
\end{equation}
However, any sequence $(0,u_1,\dots,u_n,0)$ that begins with an
increase and ends with a decrease must have exactly one more peak
than it has valleys, so the first possibility occurs $\xpe(u)$ times,
the second $\xpe(u)-1$ times, and the last $n+1-2\xpe(u)$ times.

(b) We have $C_n(t)=\aff{C}_n(1,t,\dots,t)$. The analysis is similar
to~(a), the only change being that $c_n(u)$ now specializes to
$1+t$ or $2$ according to whether $u_{n-1}<u_n$.
An equivalent way to obtain the same result would be to use the
rules in~\eqref{eq:crule} but with $u_{n+1}=\infty$.
In these terms, the sequence
$(0,u_1,\dots,u_n,\infty)$ has $\lpe(u)$ peaks each contributing
factors of $2t$, along with $\lpe(u)$ valleys each contributing
factors of $2$, and the remaining $n-2\lpe(u)$ contributions are
factors of $1+t$.

(c) A sum over $w\in C_n$ may be viewed as $2^n$ sums over
permutations of $n$ distinct objects (first choose which subset of  
letters in $[n]$ to negate). In this way, it is not
hard to see that $\aff{C}_n(1,t\dots,t,1)=2^nA_{n-1}(t)$.
The analysis of this case is similar to~(b), but now using the
convention that $u_0=u_{n+1}=\infty$.
\end{proof}

The following result was first obtained by Fulman (using the
combinatorics of shuffling~\cite{Fulman}) and later by Petersen
(using a variation of the theory of $P$-partitions~\cite{Petersen}).

\begin{cor}\label{cor:cycC}
For $n\ge1$, we have
\[
\aff{C}_n(t) = 2^n tA_{n-1}(t).
\]
\end{cor}

\begin{proof}[Proof \#1]
Comparing parts (a) and (c) of Corollary~\ref{cor:gamC},
it suffices to show that $\xpe(u)-1$ and $\pk(u)$ have the
same distribution as $u$ varies over $S_n$. To see this, recall
from the proof of Corollary~\ref{cor:gamC}(a) that every $u\in S_n$
has exactly $\xpe(u)-1$ valleys, and that these valleys occur in
internal positions. Thus $\xpe(u)-1=\pk(v)$, where $v_i=n+1-u_i$.
\end{proof}

\begin{proof}[Proof \#2]
It suffices to show that $\ades(w)$ and $\des(u)+1$ have the same
probability distributions as $w$ and $u$ vary uniformly over $C_n$
and $S_n$, respectively. To see this, first consider
\[
f_{ij}(t):=\sum_{u\in S_n:u_i=j}t^{\ades(u)},
\]
so that $\aff{A}_{n-1}(t)=f_{1j}(t)+\cdots+f_{nj}(t)$.
Since $\ades(u)$ is invariant under cyclic shifts
(recall the proof of Proposition~\ref{prp:affA}),
it follows that $f_{1j}(t)=\cdots=f_{nj}(t)$, and hence
$f_{ij}(t)=tA_{n-2}(t)$ (Corollary~\ref{cor:cycA}).
Thus, the distribution of $\ades(u)$ as $u\in S_n$ varies over
all $(n-1)!$ permutations with a fixed value in one position is
the same as the distribution of $\des(v)+1$ over $v\in S_{n-1}$.

Now consider $w\in C_n$. If we fix in advance the set $\{w_1,\ldots,w_n\}$
(one of $2^n$ equally likely possibilities), one may view the word
$\hat w:=w_1\cdots w_n0$ as a permutation of $n+1$ objects, and thus
identify $\hat w$ as one of the $n!$ members of $A_n=S_{n+1}$ with a
fixed value in its last position. However, it is not hard to see that
$\hat w$ (as a member of $A_n$) and $w$ (as a member of $C_n$) have
the same number of affine descents.
\end{proof}

\subsection{Type $B$}\label{sec:Bgam}
The Weyl group $B_n$ is identical to $C_n$, but has a root system
that is a rescaling of the $C_n$ root system; namely,
\[
\{\pm\eps_i:1\le i\le n\}\cup
  \{\pm\eps_i\pm\eps_j:1\le j<i\le n\}.
\]
We can likewise rescale the choice of simple roots; the only
change is that $\alpha_1$ is now $\eps_1$. In this way,
the positive roots are (positive) rescalings of the positive roots
for $C_n$, so ordinary descents in $B_n$ and $C_n$ are the same.
On the other hand, the lowest root $\alpha_0=-\hr$ is
now $-\eps_{n-1}-\eps_n$, so 0 is in the affine descent set
of $w\in B_n$ if and only if $w_{n-1}+w_n>0$.

For example, if $w=23\bar{4}5\bar{1}$, then $\aDes(w)=\{0,3,5\}$.

\begin{prp}\label{prp:flagB}
For $n\ge 2$, we have
\[
\aff{B}_n(t_0,\ldots,t_n)=\sum_{u \in S_n}b_1(u)\cdots b_n(u),
\]
where $b_i(u)=c_i(u)$ for $i<n-1$ as defined in~\eqref{eq:csum}, and
\begin{align*}
b_{n-1}(u)&:=t_{n-1}^{\chi(u_{n-2} < u_{n-1})}
    + (t_0t_n)^{\chi(u_{n-1} > u_n)},\\
b_n(u)&:=t_{n\vph}^{\chi(u_{n-1} < u_n)} + t_0^{\chi(u_{n-1}< u_n)}.
\end{align*}
\end{prp}

\begin{proof}
By factoring elements of $B_n$ in the form $w=\sigma u$
($\sigma\in\Z_2^n$, $u\in S_n$), the analysis proceeds as in the
proof of Proposition~\ref{prp:flagC}, except that
\begin{itemize}
\item if $u_{n-1}<u_n$, then $0\in\aDes(w)\Leftrightarrow\sigma_n=+1$,
\item if $u_{n-1}>u_n$, then $0\in\aDes(w)\Leftrightarrow\sigma_{n-1}=+1$.
\end{itemize}
Thus for each $i\in[0,n]$, there is still a unique $j$ (depending on $u$)
such that the presence or absence of $i$ in $\aDes(w)$ is controlled
by the value of $\sigma_j$. Hence,
\[
\sum_{\sigma\in\Z_2^n\vphantom{\aDes}}\ \prod_{i\in\aDes(\sigma u)}
  \!\!t_i\ = \ \left(b^-_1(u)+b^+_1(u)\right)
  \cdots \left(b^-_n(u)+b_n^+(u)\right),
\]
where $b^-_j(u)$ denotes the product of all $t_i$ such that
$\sigma_j=-1$ forces $i\in\aDes(\sigma u)$, and $b^+_j(u)$ denotes
the analogous product when $\sigma_j=+1$. These products are the
same as their counterparts for $C_n$ in Proposition~\ref{prp:flagC}
except for those involving $t_0$; namely, $b^+_{n-1}(u)$ and $b^+_n(u)$.
In the latter, there should be an extra factor of $t_0$ only when
$u_{n-1}<u_n$, and in the former there should be an extra factor
of $t_0$ when $u_{n-1}>u_n$.
\end{proof}

\newpage
\begin{cor}\label{cor:gamB}
For $n\ge 2$, we have
\[
\aff{B}_n(t) = \sum_{u \in S_n}
   \phi(u)(4t)^{\xpe(u)}(1+t)^{n+1-2\xpe(u)},
\]
where
\[
\phi(u):=
\begin{cases}
  \cbox{18}{1}   & \mbox{if } u_{n-2} > u_{n-1} > u_{n}, \\
  \cbox{18}{0}   & \mbox{if } u_{n-2} > u_{n} > u_{n-1}, \\
  \cbox{18}{1/2} & \mbox{otherwise.}
\end{cases}
\]
In particular, $\aff{B}_n(t)$ has a nonnegative $\gamma$-vector.
\end{cor}

Given our convention that $u_0=0$, one should understand that
$\phi(u)=1/2$ for $u\in S_2$.

\begin{proof}
Recall from the proof of Corollary~\ref{cor:gamC}(a) that if we
specialize the variables so that $t_i\to t$ for all $i$, we obtain
\begin{equation}\label{eq:cspec}
  c_1(u)\cdots c_n(u)\to (1/2)(4t)^{\xpe(u)}(1+t)^{n+1-2\xpe(u)}.
\end{equation}
Comparing the definitions of $b_i(u)$ and $c_i(u)$, we see that
\[
\frac{b_1(u)\cdots b_n(u)}{c_1(u)\cdots c_n(u)}
  =\frac{b_{n-1}(u)b_n(u)}{c_{n-1}(u)c_n(u)}\to
\begin{cases}
   2(1+t^2)/(1+t)^2&\text{if $u_{n-2}>u_{n-1}>u_n$},\\
   \cbox{79}{1}&\text{otherwise}.
\end{cases}
\]
Thus, $b_1(u)\cdots b_n(u)$ usually specializes in the same way
as in~\eqref{eq:cspec}.

Now consider that transposing $u_{n-1}$ and $u_n$ yields
a bijection $u\leftrightarrow u'$ between permutations in $S_n$
that satisfy $u_{n-2}>u_{n-1}>u_n$ and $u'_{n-2}>u'_n>u'_{n-1}$.
Furthermore, we have
\[
\frac{b_1(u')\cdots b_n(u')}{c_1(u)\cdots c_n(u)}
  =\frac{b_{n-1}(u')b_n(u')}{c_{n-1}(u)c_n(u)}
  \to\frac{4t}{(1+t)^2}.
\]
Therefore, if we combine the terms indexed by $u'$ and $u$ in these
cases (and eliminate $u'$ from the sum), the net contribution
of $u$ is $(2(1+t^2)+4t)/(1+t)^2=2$ times~\eqref{eq:cspec}.
\end{proof}

\subsection{Type $D$}\label{sec:Dgam}
The Weyl group $D_n$ is the subgroup of $B_n$ consisting
of signed permutations $w=w_1\cdots w_n$ with an even number of
negative entries. It has a root system of the form
\[
\{\pm\eps_i\pm\eps_j:1\le j<i\le n\},
\]
and one can choose simple roots so that $\alpha_1=\eps_2+\eps_1$
and $\alpha_i=\eps_i-\eps_{i-1}$ for $2\le i\le n$. This choice
is compatible with our previous choices for $B_n$ and $C_n$ in
the sense that a $D_n$ root is positive if and only if it is
positive as a $B_n$ or $C_n$ root.

It is important to note that $D_n$ is irreducible only for $n\ge3$.
In such cases, the lowest root $\alpha_0=-\hr$ is $-\eps_{n-1}-\eps_n$
(the same as in $B_n$), and thus the affine descent set of $w\in D_n$
consists of all $i\in[2,n]$ such that $w_{i-1}>w_i$, together with 1
(if $w_1+w_2<0$) and $0$ (if $w_{n-1}+w_n>0$).

For example, if $w=3\bar{4}2\bar{1}5$, then $\aDes(w)=\{0,1,2,4\}$.

Note that by specializing $t_0=1$ in the following, we recover
Stembridge's formula for the flag $h$-polynomial of the Coxeter
complex $\Sigma(D_n)$ (Proposition~A.4 in~\cite{Stembridge}).

\begin{prp}\label{prp:flagD}
For $n\ge 4$, we have
\[
\aff{D}_n(t_0,\ldots,t_n) = \frac{1}{2}
   \sum_{u \in S_n} d_1(u)\cdots d_n(u),
\]
where $d_i(u)=b_i(u)$ for $i>2$ as defined in
Proposition~\ref{prp:flagB}, and
\begin{align*}
d_1(u)&:=t_1^{\chi(u_1>u_2)} + t_2^{\chi(u_1 > u_2)},\\
d_2(u)&:=(t_1t_2)^{\chi(u_1<u_2)} + t_3^{\chi(u_2>u_3)}.
\end{align*}
\end{prp}

\begin{proof}
Following the proof of Proposition~A.4 in~\cite{Stembridge}, note that
the definition of an affine descent set in $D_n$ makes sense for any
signed permutation $w\in B_n$. Since replacing $1$ with $-1$ or
vice-versa in $w_1\cdots w_n$ does not change this set, it follows that
\[
\aff{D}_n(t_0,\dots,t_n)=\frac{1}{2}\sum_{u\in S_n\vphantom{\aDes}}
  \ \sum_{\sigma\in\Z_2^n\vphantom{\aDes}}\ \prod_{i\in\aDes(w)}t_i.
\]
The analysis of $w=\sigma u$ now proceeds as in the proof of
Proposition~\ref{prp:flagB}, except that
\begin{itemize}
\item if $u_1>u_2$, then $1\in\aDes(w)\Leftrightarrow\sigma_1=-1$,
\item if $u_1<u_2$, then $1\in\aDes(w)\Leftrightarrow\sigma_2=-1$.
\end{itemize}
Again it follows that for each $i\in[0,n]$, there is a unique $j$
(depending on $u$) such that the presence or absence of $i$ in
$\aDes(w)$ is controlled by the value of $\sigma_j$. Hence,
\[
\sum_{\sigma\in\Z_2^n\vphantom{\aDes}}\ \prod_{i\in\aDes(w)}t_i
 \ =\ \left(d^-_1(u)+d^+_1(u)\right)\cdots\left(d^-_n(u)+d^+_n(u)\right),
\]
where $d_j^-(u)$ denotes the product of all $t_i$ such that
$\sigma_j=-1$ forces $i\in\aDes(\sigma u)$, and $d_j^+(u)$ denotes
the analogous product when $\sigma_j=+1$. These products are the same
as their counterparts for $B_n$ in Proposition~\ref{prp:flagB}
except for those that involve $t_1$; namely, $d_1^-(u)$ and $d_2^-(u)$.
In the former, there should be a factor of $t_1$ only when $u_1>u_2$,
and in the latter, there should be an extra factor
of $t_1$ when $u_1<u_2$.
\end{proof}

Specializing, we obtain nonnegative $\gamma$-expansions for
both $\aff{D}_n(t)$ and $D_n(t)$, the latter of which is due to
Stembridge (Corollary~A.5 in~\cite{Stembridge};
compare also Theorem~6.9 in~\cite{Chow}).

\begin{cor}\label{cor:gamD}
For $n \ge 4$, we have
\begin{itemize}\vspace{2pt}
\item[(a)] $\aff{D}_n(t) = \sum_{u \in S_n} \phi(u)
   \phi(\rev{u})(4t)^{\xpe(u)}(1+t)^{n+1-2\xpe(u)}$.\vspace{6pt}
\item[(b)] $D_n(t)=\sum_{u \in S_n}
   \phi(\rev{u})(4t)^{\lpe(u)}(1+t)^{n-2\lpe(u)}$.\vspace{3pt}
\end{itemize}
where $\phi(u)$ is defined as in Corollary \ref{cor:gamB} and
$\rev{u}:= u_n \cdots u_2u_1$.
\end{cor}

\begin{proof}
(a) Specializing the variables so that $t_i\to t$ for all $i$,
we obtain
\[
\frac{d_1(u)\cdots d_n(u)}{b_1(u)\cdots b_n(u)}
  =\frac{d_1(u)d_2(u)}{b_1(u)b_2(u)}\to
\begin{cases}
   2(1+t^2)/(1+t)^2&\text{if $u_1<u_2<u_3$},\\
   \cbox{79}{1}&\text{otherwise}.
\end{cases}
\]
Now pair each permutation $u\in S_n$ such that $u_1<u_2<u_3$
with the permutation $u''$ obtained by switching $u_1$ and $u_2$.
In such cases, we have
\[
\frac{d_1(u'')\cdots d_n(u'')}{b_1(u)\cdots b_n(u)}
  =\frac{d_1(u'')d_2(u'')}{b_1(u)b_2(u)}
  \to\frac{4t}{(1+t)^2},
\]
so when the expansion in Proposition~\ref{prp:flagD} is specialized,
the terms indexed by $u$ and $u''$ such that $u_1<u_2<u_3$ may be
combined into a single term with twice the $b$-weight of $u$, yielding
\[
\sum_{u\in S_n}\phi(\rev{u})b_1(u)\cdots b_n(u)\to \aff{D}_n(t).
\]
Now proceed as in the proof of Corollary~\ref{cor:gamB}, combining the
terms indexed by $u\in S_n$ such that $u_{n-2}>u_{n-1}>u_n$ with the
terms indexed by the permutations $u'$ obtained by switching $u_{n-1}$
and $u_n$, and note that $\phi(\rev{u})=\phi(\rev{(u')})$, even when $n=4$.

(b) Similarly, we have $D_n(t)=\aff{D}_n(1,t,\dots,t)$.
Under this specialization, the effects on $b_i(u)$ and $d_i(u)$ are
similar to the previous case; the only differences occur in the terms
that involve $t_0$; namely $b_i(u)$ and $d_i(u)$ for $i=n-1$ and $i=n$.
However, we have $b_i(u)=d_i(u)$ in these cases (even without
specialization), so the same reasoning as above implies
\[
\sum_{u\in S_n}\phi(\rev{u})b_1(u)\cdots b_n(u)
  \Big\vert_{t_0=1}\to \ D_n(t).
\]
Now observe that $b_i(u)=c_i(u)$ for all $i$ when $t_0=1$,
so the result follows by the reasoning in the proof of
Corollary~\ref{cor:gamC}(b).
\end{proof}

\section{Identities and generating functions}\label{sec:ident}
\subsection{Strange Identities}
Here we provide several unexpected identities (two new, one old)
relating the ordinary and affine Eulerian polynomials.

\begin{prp}\label{prp:CBCid}
For $n\ge2$, we have
\[
2\aff{C}_n(t) = \aff{B}_n(t) + 2nt C_{n-1}(t).
\]
\end{prp}

\begin{proof}
Given $u\in S_n$, let $\trim{u}=u_1\cdots u_{n-1}$, a permutation
of $n-1$ distinct positive integers. Noting that the definitions
of peak numbers make sense for any sequence of positive integers,
we see that the distribution of $\lpe(\trim{u})$ as $u$ varies
over $S_n$ is the same as $n$ copies of the distribution of $\lpe(v)$
as $v$ varies over $S_{n-1}$. Thus Corollary~\ref{cor:gamC}(b) implies
\[
2ntC_{n-1}(t)=2t \sum_{u \in S_n}
  (4t)^{\lpe(\trim{u})}(1+t)^{n-1-2\lpe(\trim{u})}.
\]
Now recall that swapping $u_{n-1}$ and $u_n$ provides a bijection
between the permutations $u\in S_n$ satisfying $\phi(u)=1$
(i.e., $u_{n-2}>u_{n-1}>u_n$) with the permutations $u'$
satisfying $\phi(u')=0$ (i.e., $u'_{n-2}>u'_n>u'_{n-1}$).
Noting that $\lpe(\trim{u})=\lpe(\trim{u'})$ for such pairs,
we can achieve an equivalent result by doubling the contribution
of $u'$ and eliminating $u$, or simply modify the contribution
of every permutation $u$ by the factor $2(1-\phi(u))$. Thus,
\begin{equation}\label{eq:Cid}
2ntC_{n-1}(t)=\sum_{u \in S_n}(1-\phi(u))
  (4t)^{\lpe(\trim{u})+1}(1+t)^{n-1-2\lpe(\trim{u})}.
\end{equation}
On the other hand, Corollaries~\ref{cor:gamC}(a) and \ref{cor:gamB} imply
\[
2\aff{C}_n(t)-\aff{B}_n(t) = \sum_{u \in S_n}
  (1-\phi(u))(4t)^{\xpe(u)}(1+t)^{n+1-2\xpe(u)}.
\]
Noting that $\xpe(u)=\lpe(\trim{u})+1$ whenever $\phi(u)\ne1$,
the result follows.
\end{proof}

\begin{prp}\label{prp:BDDid}
For $n\ge3$, we have
\[
\aff{B}_n(t) = \aff{D}_n(t) + 2nt D_{n-1}(t).
\]
\end{prp}

\begin{proof}
It is easy to check the case $n=3$ (note that $D_2(t)=(1+t)^2$),
so we assume $n\ge4$.

Following the argument we used to deduce~\eqref{eq:Cid}
from Corollary~\ref{cor:gamC}(b), one may similarly use
Corollary~\ref{cor:gamD}(b) to show that
\begin{align*}
2ntD_{n-1}(t)
  &=\sum_{u \in S_n}(1-\phi(u))\phi(\rev{u})
    (4t)^{\lpe(\trim{u})+1}(1+t)^{n-1-2\lpe(\trim{u})}\\
  &=\sum_{u \in S_n}(1-\phi(u))\phi(\rev{u})
    (4t)^{\xpe(u)}(1+t)^{n+1-2\xpe(u)},
\end{align*}
again using the fact that $\xpe(u)=\lpe(\trim{u})+1$
when $\phi(u)\ne1$. The only caveats are that one needs to
check that $\phi(\rev{(\trim{u})})=\phi(\rev{u})$ for all $u\in S_n$,
and $\phi(\rev{u})=\phi(\rev{(u')})$ when $\phi(u)=1$.
One should also check that the formula provided in
Corollary~\ref{cor:gamD}(b) is valid for $D_3$,
since the argument given there is not.

On the other hand,
Corollaries~\ref{cor:gamB} and \ref{cor:gamD}(a) imply
\[
\aff{B}_n(t)-\aff{D}_n(t) = \sum_{u \in S_n}
  (1-\phi(\rev{u}))\phi(u)(4t)^{\xpe(u)}(1+t)^{n+1-2\xpe(u)}.
\]
Comparing the two expansions and noting that
$\xpe(\rev{u})=\xpe(u)$, the result follows.
\end{proof}

The following identity is due to Stembridge
(set $l=0$ in \cite[Lemma~9.1]{Stembridge2}).

\begin{prp}\label{prp:BDAid}
For $n\ge 2$, we have
\[
B_n(t)=C_n(t) = D_n(t) + n2^{n-1}tA_{n-2}(t).
\]
\end{prp}

\begin{proof}
It is easy to check the cases $n=2$ and $n=3$, so assume $n\ge4$.

By the same reasoning we used in the proof of
Proposition~\ref{prp:CBCid}, Corollary~\ref{cor:gamC}(c) implies
\[
n2^{n-1}tA_{n-2}(t) = 2t \sum_{u \in S_n}
  (4t)^{\pk(\ltrim{u})}(1+t)^{n-2-2\pk(\ltrim{u})},
\]
where $\ltrim{u}:=u_2\cdots u_n$. Now consider that
if $\phi(\rev{u})=1$ (i.e., $u_1<u_2<u_3$) and
$u''$ is obtained from $u$ by switching $u_1$ and $u_2$
(hence $\phi(\rev{(u'')})=0$) then $\pk(\ltrim{u})=\pk(\ltrim{u''})$.
It follows that we can create an equivalent sum by modifying
the contribution of every permutation $u$ by the factor
$2(1-\phi(\rev{u}))$, yielding
\begin{align*}
n2^{n-1}tA_{n-2}(t)
  &= \sum_{u \in S_n}(1-\phi(\rev{u}))
    (4t)^{\pk(\ltrim{u})+1}(1+t)^{n-2-2\pk(\ltrim{u})}\\
  &= \sum_{u \in S_n}(1-\phi(\rev{u}))
    (4t)^{\lpe(u)}(1+t)^{n-2\lpe(u)},
\end{align*}
using the fact that $\lpe(u)=\pk(\ltrim{u})+1$
whenever $\phi(\rev{u})\ne1$. On the other hand,
it is clear from Corollaries~\ref{cor:gamC}(b)
and \ref{cor:gamD}(b) that this sum is $C_n(t)-D_n(t)$.
\end{proof}

\subsection{Generating functions}
First let us review the (known) generating functions for the Eulerian
polynomials corresponding to the Weyl groups $A_n$, $B_n=C_n$, and $D_n$:
\begin{align}
A(t,z)&:=1+\sum_{n \ge 1} tA_{n-1}(t) \frac{z^n}{n!}
  = \frac{(1-t)}{1-te^{z(1-t)}},\label{eq:Agf}\\
B(t,z)=C(t,z)&:=  1+ (1+t)z + \sum_{n \ge 2 } B_n(t) \frac{z^n}{n!}
  = \frac{(1-t)e^{z(1-t)} }{ 1-te^{2z(1-t)} },\label{eq:BCgf}\\
D(t,z)&:=  1+ tz + \sum_{n \ge 2 } D_n(t) \frac{z^n}{n!}
  = \frac{(1-t)(e^{z(1-t)}-z)}{1-te^{2z(1-t)}}.\label{eq:Dgf}
\end{align}
The first of these is classical (e.g., see~\cite[p.~244]{Comtet}),
and for proofs of the latter two see Theorem~3.4 and Corollary~4.9
in~\cite{Brenti}, but note that our initial terms for $D(t,z)$, and
hence the resulting closed form, are slightly different from those
in~\cite{Brenti}.

For small $n$, the values for $A_n(t)$, $B_n(t)$, $D_n(t)$
implicit in the above expansions are
\[
A_{-1}(t)=1/t, \ B_0(t)=1,\ B_1(t)=1+t,\ D_0(t)=1, \ D_1(t)=t.
\]
In this way, Proposition~\ref{prp:BDAid} is valid even for $n=0$
or $1$, and immediately implies
\[
B(t,z)=D(t,z)+zA(t,2z).
\]
Thus~\eqref{eq:Dgf} may be viewed as a corollary of~\eqref{eq:BCgf}
and Proposition~\ref{prp:BDAid}.

Turning to the affine Eulerian polynomials, let
\begin{align*}
\aff{A}(t,z)&:= z+\sum_{n\ge 2} \aff{A}_{n-1}(t) \frac{z^n}{n!},\\
\aff{C}(t,z)&:= 1+\sum_{n\ge1}\aff{C}_n(t)\frac{z^n}{n!},\\
\aff{B}(t,z)&:=  2 + 2tz + \sum_{n \ge 2 }
  \aff{B}_n(t) \frac{z^n}{n!},\\
\aff{D}(t,z)&:=  2 + 4t\frac{z^2}{2} + \sum_{n \ge 3 }
  \aff{D}_n(t) \frac{z^n}{n!}.
\end{align*}

\begin{prp}\label{prp:aBDgf}
We have
\begin{align}
\aff{A}(t,z)&= \frac{ z(1-t) }{1-te^{z(1-t)}},\\
\aff{C}(t,z)&= \frac{1-t}{1-te^{2z(1-t)}},\label{eq:aCgf}\\
\aff{B}(t,z)&=\frac{2(1-t)(1-tze^{z(1-t)}) }{ 1-te^{2z(1-t)} },\\
\aff{D}(t,z)&=\frac{ 2(1-t)(1+tz^2-2tze^{z(1-t)})}{ 1-te^{2z(1-t)} }.
\end{align}
\end{prp}

\begin{proof}
Corollaries~\ref{cor:cycA} and~\ref{cor:cycC} immediately
imply $\aff{A}(t,z)=zA(t,z)$ and $\aff{C}(t,z)=A(t,2z)$,
so these generating functions are consequences of~\eqref{eq:Agf}.
In the remaining two cases,
the values implicit for small $n$ in the series defined above
(namely, $\aff{B}_0(t)=\aff{D}_0(t)=2$, $\aff{C}_0(t)=1$,
$\aff{B}_1(t)=\aff{C}_1(t)=2t$, $\aff{D}_1(t)=0$, and $\aff{D}_2(t)=4t$)
have been deliberately chosen so that Propositions~\ref{prp:CBCid}
and~\ref{prp:BDDid} remain valid for all $n\ge0$,
and thus respectively imply
\begin{align*}
2\aff{C}(t,z)=\aff{B}(t,z)+2tzC(t,z),\\
\aff{B}(t,z)=\aff{D}(t,z)+2tzD(t,z).
\end{align*}
The first of these, together with~\eqref{eq:BCgf} and~\eqref{eq:aCgf},
yields the formula claimed for $\aff{B}(t,z)$, and then the second
(with~\eqref{eq:Dgf}) yields the formula claimed for $\aff{D}(t,z)$.
\end{proof}

}

\end{document}